\documentclass[11pt]{article}

\usepackage{amsfonts}
\usepackage{amscd}
\usepackage{amsthm}
\usepackage{indentfirst}
\usepackage[margin=3cm]{geometry}
\usepackage[small,nohug]{diagrams}
\diagramstyle[labelstyle=\scriptstyle]

\usepackage[all]{xy}
\usepackage{url}

\usepackage{amssymb, amsmath, amscd, amsthm, amsgen, amstext, amsbsy, amsopn, diagrams}
\usepackage{dsfont}
\usepackage{epsfig}
\usepackage{epic,eepic}

\usepackage{verbatim}

\theoremstyle{definition}

\newcommand{\G}[1][\R^{2n} \times S^1] {{\text{\rm Cont}_0(#1)}}
\newcommand{\Ham}[1][\R^{2n}]{{\text{\rm Ham}(#1)}}

\newcommand{\GZK}[1][\R^{2n} \times S^1] {{\text{\rm Cont}_0^{\Z_k}(#1)}}

\newcommand{\R}{{\mathbb{R}}}

\newcommand{\Z}{{\mathbb{Z}}}
\newcommand{\N}{{\mathbb{N}}}

\begin{document}

\author{\textsc Maia Fraser\thanks{Department of Mathematics, University of Toronto} }

\title{Contact spectral invariants and persistence}

\date{January 28, 2014}
\maketitle

\begin{abstract}
This sketch shows that the usual generating function based capacities have an interpretation in 
the language of persistent homology as persistences of certain homology classes in 
the persistence module formed by the corresponding generating function 
homology groups. This viewpoint suggests various new invariants, in particular
a $\Z_k$-equivariant capacity which can be used to prove orderability of
lens spaces, proved by Milin \cite{milin} using contact homology 
and by Sandon \cite{sandon-equiv}
using equivariant generating function homology. These are informal notes 
originally circulated in January 2014.
\end{abstract}

\section{Generating-function based capacities} \label{sec:cap-pers}

The usual Viterbo/Sandon capacity \cite{viterbo, sandon} is defined for $\phi \in \Ham$,
resp. $\G$, by first assuming a
quadratic at infinity generating function $S: E \to \R$ for the
(exact Lagrangian, resp. Legendrian) ``graph" of $\phi$.
Here $E$ is a vector bundle over $M =
S^{2n}$ resp. $S^{2n} \times S^1$, the original manifold compactified,
and may even be taken to be trivial.
One then considers sub-level sets $E^a :=\{S \leq a \} \subseteq E$ and puts
$$c(\phi) := \inf \{a \in \R : \theta_a(\mu) \neq 0 \}$$
where $\mu \in H^{2n}(M)$ is the orientation class of the base
$M$
and for any $a \in \R$, $\theta_a: H^*(M) \to H^*(E^a, E^{-})$ is 
the composition 
$\theta_a = i^*_a \circ\, \theta$ 
of the natural isomorphism 
(c.f. Thom isomorphism
and Excision theorem)\footnote
{
Technically, the degree $*$ is shifted by $j$, the index of $Q_\infty$, the quadratic at infinity part of $S$.
Indeed by the Thom isomorphism, $H^*(M) \cong H^{*-j}(DE_{-}, SE_{-})$, 
where $E_{-}$ is any rank $j$ real vector bundle over $M$.
This is isomorphic to
$H^{*-j}(E_{-}, CE_{-})$, letting $CE_{-}$ be the closure
of $(DE_{-})^c \subset E_{-}$, and excising $(DE_{-})^c$.
On the other hand, taking $E_{-}$ to be  the
the maximal sub-bundle on which $Q_\infty$ is negative definite,
$H^{*-j}(E_{-}, CE_{-}) \cong H^{*-j}(E, E^{-\infty})$ because
for sufficiently negative $b$, $E^{b}$ and $E^{-\infty}$ are homotopy equivalent 
and the pair $(E, E^b)$ fiberwise deformation retracts to the pair
$(E_-, \{e \in E_{-} : |e| \geq B\})$ for some $B \in \R$.
}
$\theta : H^*(M) \to H^*(E, E^{-\infty})$ 
with the homomorphism
$i^*_a : H^*(E, E^{-\infty}) \to H^*(E^a, E^{-\infty})$ induced by inclusion
$i_a: E^a \to E$.

\medskip
{\bf Assume coefficients in a field.}
Note: 
\begin{itemize}
\item $c(\phi)$, resp. $[c(\phi)]$ in the contact case\footnote{Here $[x]$ denotes the largest $a \in \Z$ s.t. $a < x$; any $c(\phi, \mu)$, resp. $[c(\phi, \mu)]$ is also conjugation-invariant.}, is conjugation-invariant
\item An analogous $c(\phi) = c(\phi, \mu)$ can be defined using any non-trivial $\mu \in H^*(M)$. This does not give any more information\footnote{
Indeed $H^*(S^{2n})$ is non-trivial only for $*=0, 2n$
and cyclic in these cases (while the generator of $H^0(S^{2n})$ is killed
at $a=0$ for all compactly supported $\phi$).
} about $\phi$
in the symplectic case. Nor does it in the contact case for $\phi \in \G$ which is a lift  $\phi= \bar \sigma$ of 
$\sigma \in \Ham$ because, as may be checked, $c_\text{contact}(\bar\sigma) = [c_\text{symp.}(\sigma)]$.
\item By letting $\eta = \theta(\mu) \in H^p(E, E^{-\infty})$ ($p = 2n$ resp. $2n+1$) we may equivalently write
\begin{align*}
c(\phi) &:= \inf \{a \in \R : i^*_a(\eta) \neq 0 \}.
\end{align*}
Alternatively, 
\begin{align}
c(\phi) &:= \sup \{a \in \R : (j_a)_*(\eta) \neq 0 \} \label{eq:eta}
\end{align}
where $(j_a)_* : H_p(E, E^{-\infty}) \to H_p(E, E^a)$ induced by inclusion
$j_a: E^{-\infty} \to E^a$ with $\eta = \theta([M])$, $\theta$ re-defined analogously for homology, and
$[M]$ being the fundamental class of the base, i.e. Kronecker dual to $\mu$.
\item This is in fact the {\it persistence} of the class $\eta$, in the sense of persistent homology (see next Section), for the {\it persistence module} 
\begin{align*}
(V_{a} ~\vert~ a \in \R), \qquad V_{a} := H_p(E, E^{a}) 
\end{align*}
(with linear map $V_a \to V_b$ induced by inclusion whenever $a \leq b$) and moreover 
many more such persistences and persistent homology groups can be defined; under certain circumstances they too will be conjugation-invariant.
\end{itemize}

\section{Persistence}
A {\bf persistence module} $\mathbb V$ over $\R$ is 
an indexed family of vector spaces
$$(V_t ~\vert~ t \in \R)$$
and a doubly-indexed family of linear maps
$$(v^t_s \colon V_s \to V_t ~\vert~ s \leq t)$$
which satisfy the composition law
$$v^t_s \circ v^s_r = v^t_r$$
whenever $r \leq s \leq t$ and where $v^t_t$ is the identity map on $V_t$. 
For more on persistent homology see Weinberger \cite{shmuel}. For the framework of persistent modules see Chazal et al \cite{vin} who introduced this useful formalism.
Equivalently, a persistence module can be viewed as a functor from $\R$ (as a category 
with a unique morphism $s \to t$ whenever $s \leq t$) to the category of vector spaces.
Let us call
the maps $v^t_s$ {\bf persistence maps} (not standard terminology).

Any functor from topological spaces to vector spaces, for example singular 
homology in a fixed degree
$H(\cdot) := H_p(X, K), ~K$ a field, produces a persistence module when applied to 
a topological space with an $\R$-filtration, $X = \cup_{t \in \R} X^t,
s \leq t \Rightarrow X^s \subseteq X^t$ using as linear maps $v^t_s$ the homomorphisms induced
by inclusion. A contravariant functor, such as singular cohomology could be used instead; this
produces a persistence module if the indexing of
the filtration is first reversed (so $s \leq t \Rightarrow X^s \supseteq X^t$).

\smallskip
{\sc Remark:} In principle, equivariant (co)homology groups could also be used (when sets $X^t$
of the filtration are invariant under an ambient group action). 
To my knowledge this has not been done in topological data applications.
An application in contact/symplectic geometry is proposed in
Section 4.

\smallskip 
Note that defining a persistence module via singular homology with coefficients in a field
also produces an intermediate persistence module, singular chain
groups $C_p(X^t), ~t \in \R$ which
form an {\bf $\R$-family of chain complexes} $C_\bullet (X^t), ~t \in \R$
where the chain groups are vector spaces
and
for any $s, t \in \R$ the persistence maps $v^t_s := (i_s^t)_*$ (in all degrees, induced by
inclusion $i_s^t: X_s \to X^t$) define a chain map
between the complexes at level $s$ and $t$.
\begin{diagram}
\ldots C_{p+1}(X^s) & \rTo^{d_{p+1}} & C_p(X^s) & \rTo^{d_p} & C_{p-1}(X^s) \ldots \\
\dTo^{v^t_s} &&\dTo^{v^t_s} && \dTo^{v^t_s} \\
\ldots C_{p+1}(X^t) &  \rTo^{d_{p+1}} &C_p(X^t) &  \rTo^{d_p} & C_{p-1}(X^t) \ldots
\end{diagram}
Such an object can also
be used as starting point. 
Taking homology one then obtains persistence modules $H_p(X^t), ~t \in \R$ 
whose union over all $p \in \Z$
forms what can be thought of as a $\Z$-graded persistence module (non-standard terminology). 

\smallskip
{\sc Remark:} Suppose there
is an abelian group $G$ acting on each of the chain complexes in the
{$\R$-family of chain complexes} $C_\bullet (X^t), ~t \in \R$; i.e., for each $t \in \R$, 
and each $p \in \Z$, $G$ acts on the vector space $C_p (X^t)$ and the differential
$d_{p+1}:C_{p+1}(X^t) \to C_ p(X^t)$ is $G$-equivariant. Then we
can define\footnote{Consider $C_\bullet (X^t)$ as a chain complex of $K[G]$-modules, where $K[G]$ denotes
the group ring of $G$ (in our case the group algebra, since $K$ is a field). 
$G$-equivariant homology of $C_\bullet (X^t)$ is computed
by taking a projective resolution $(E_\bullet, \delta)$ of $K$ as a $K[G]$-module (with trivial
$G$-action), tensoring 
$(E_\bullet, \delta)$ with $(C_\bullet, d)$ and taking the homology of 
$E_\bullet \otimes_{K[G]} C_\bullet$ (with suitable differential).} 
$G$-equivariant homology groups $H_{G, p}(X^t), ~t \in \R$. 
If, in addition to the $G$-action on each chain complex $C_\bullet (X^t)$ for $t \in \R$ 
we require that
persistence maps $v_s^t$ define morphisms in the category of chain complexes with
$G$-actions, i.e. commute not only with the differential but also the group action,
then the $G$-equivariant homology groups $H_{G, p}(X^t), ~t \in \R$
will form a persistence module in each degree $p\in \Z$. 
This construction has not been considered in topological data analysis but a
possible application in computational geometry is work in progress.
\smallskip

{\sc Terminology:}
Consider 
a persistence module defined by homology in degree $p$.
The image $v^s_{s+a}(H_p(X^s)) \subseteq H_p(X^{s+a})$ for $a>0$ 
is called the {\bf $a$-persistent homology group} of $X^s$ in degree $k$.
An analogous object for a general persistence module would also make sense
(though no general name exists).
Further terminology which does apply to general
persistence modules, but which we give here in terms of homology for simplicity
is as follows.
For any non-trivial class $\mu \in H_p(X^s)$ its {\bf persistence}
is $$\rho(\mu) := \sup \{a >0 : v_s^{s+a}(\mu) \neq 0 \}.$$
More commonly, this information is collected for all $\mu \in H_p(X^s), s \in \R$ and recorded as $\mathbb{P}_p$, a
multiset of pairs, each pair specifying birth and death for some $\mu$:
\begin{align*}
\mathbb{P}_p := \{ (s, a+s) \in \R^2 : &(\exists \mu \in H_p(X^s)) \\
&\mu \notin v^s_r (H_p(X^r)) \text{ for }r < s, \\
&\rho(\mu) = a\}.
\end{align*}
Arranged as
a set of points in the plane with multiplicities,
$\mathbb{P}_p$ is referred to as the {\bf persistence diagram}
in degree $p$ of $X$. Alternatively this information is sometimes recorded
as a {\bf barcode}: a collection of 
horizontal bars stacked (in no particular order)
above the $x$-axis in the plane, with one bar for each $\mu$,
having left endpoint $s$ and right endpoint $s+a$ as specified above.
We think of $s$ and $s+a$ as the {birth} and {death} of $\mu$ and the length of the bar thus represents
the lifespan of $\mu$.
Stability theorems exist in various settings. Roughly speaking, they
show that when the function $f$ defining the
filtration is changed by less than $C$, classes with lifespans longer than 
$C$ will continue to exist (their endpoints will not be shifted by more than $C$). 
See \cite{shmuel} for an interesting application
of this (proving a Theorem of Gromov).

\section{GF-based capacities as persistences} \label{sec:suggestions}
\paragraph{\sc I) Extension of GF homology group definition}
Traynor \cite{traynor} defines
generating function homology groups
$G^{(a, b]}_* (\phi)$, $\phi \in \Ham$ for action windows $(a, b]$, $a < b \in \R$ {\it both nonzero}
and requires the generating function for $\phi$ not to have $a$ or $b$ as critical value.
In the contact setting, for $\phi \in \G$, 
Sandon \cite{sandon, sandon-equiv} makes the same requirement on critical values
(but does not explicitly state that $a, b$ be nonzero,
just that they be integer in order to obtain conjugation-invariance). In any case, in all these works,
GF homology groups are only computed for action windows $(a, \infty]$ 
where $a >0$.
The non-zero condition is needed since
generating functions for 
compactly supported $\phi$ always have zero as a critical value,
but $G^{(a, b]}_* (\phi)$ is defined as relative homology 
$H_*(\{S \leq b\}, \{S \leq a\})$ 
of sub-level sets of a generating function $S$ for $\phi$
having $a$ and $b$ as regular values.

For use in the constructions below, it will be convenient to 
define\footnote{An alternate definition $G_*^{(0_+, \infty]}(\phi) := G_*^{(\epsilon_\phi, \infty]}(\phi)$ for sufficiently small $\epsilon_\phi > 0$ was stated in the first draft of these notes but to do so properly requires some care with the allowed classes of $\phi$ in order to ensure generic generating functions; the current definition avoids these technicalities.} a surrogate $G_*^{(0_+, \infty]}(\phi)$
for $G^{(0, \infty]}(\phi)$. 
In fact, $G^{(-\infty, \infty]}_*(\phi) = G^{(a, \infty]}_*(\phi)$ coincides for all $\phi$
and all $a < 0$,
and is nonzero in degrees 0 and $2n$ only. We denote the generator in degree 
$0$ by $\kappa$, the generator in degree $2n$ by $\eta$, as in
Section~\ref{sec:cap-pers}. Note that
$c(\phi, \kappa) = 0$. 
We therefore quotient out
the subspace 
$\langle \kappa \rangle$ from the $\Z$-graded vector space $G^{(-\epsilon, \infty]}_*(\phi)$,
to define the surrogate, $G_*^{(0_+, \infty]}(\phi) := G^{(-\epsilon, \infty]}_*(\phi)/\langle\kappa\rangle.$
Since $c(\phi, \mu)$ is conjugation invariant for any $\mu$, the vector space
$G_*^{(0_+, \infty]}(\phi)$ is by definition conjugation-invariant. We denote
$\nu_\epsilon$ the projection to the quotient.
In the next Section we will consider another
system of GF homology groups and create a similar surrogate but in that case there
will be more generators such that $c(\phi, \mu) = 0$, and we quotient
all of them out to define the surrogate.

\paragraph{\sc II) Viterbo/Sandon capacities re-expressed}
For each $\phi$, and fixed $p \in \Z$,
let 
$$V_{a} (\phi) := \left \{
\begin{array}{ll}
G_p^{(a, \infty]}(\phi)  &\text{ if } a < 0 \\ 
G_p^{(0_+, \infty]}(\phi) &\text{ if } a=0  \\
G_p^{(a, \infty]}(\phi) &\text{ if } a>0
\end{array}
\right .
$$
Then
$(V_a(\phi) ~\vert~ a \in \R)$ is a persistence module
with persistence maps $(v_a^b): V_{a}(\phi) \to V_{b}(\phi)$, $a < b$,
induced by the inclusion
map $i_a^b \colon (E, E^a) \to (E, E^b)$ on pairs with the modification
that $v_0^b$ is defined as $\overline{(i_{-\epsilon}^b)_*}$, the linear 
map induced on the quotient by the linear map $(i_{-\epsilon}^b)_*$ with $-\epsilon \in (-\infty, 0)$ arbitrarily chosen,
and $v_a^0$ is defined as $\nu_\epsilon \circ(i_a^{-\epsilon})_*$  with $-\epsilon \in (a, 0)$
arbitrarily chosen.

\smallskip
Moreover, the Viterbo/Sandon capacity $c(\phi)$ defined
in the first Section can equivalently be expressed as the persistence
$$c(\phi) := \rho(\eta')$$
of $\eta' := v_{-\infty}^0 (\eta) \in V_0$. This is because 
$c(\phi, \eta) > 0$ for all $\phi$ \cite{traynor, sandon} (recall the definition of $\eta$ in \eqref{eq:eta}) and $v_0^b(\eta') \neq 0$ 
if and only if $(i_{-\epsilon}^b)_*(\eta') \neq 0$
which in the notation of first Section is equivalent to $(j_b)_*(\eta) \neq 0$.

\section{$\Z_k$-equivariant capacity}

In this section we discuss how to define a $\Z_k$-equivariant capacity
by using the $\Z_k$-equivariant GF homology groups of Sandon \cite{sandon-equiv}
to define a persistence module as above. 

\paragraph{\sc I)  Contact version} 
The equivariant GF homology groups 
$$
G_{\Z_k, p}^{(a, b]} (\phi) := H_{\Z_k, p}(E^a, E^b), \text{ for }a, b \in \Z 
$$
defined by Sandon in $\R^{2n} \times S^1$ \cite{sandon-equiv} for $\phi \in \GZK$
(the identity component of the group of $\Z_k$-equivariant
contactomorphisms with compact support) 
are invariant under conjugation by $\psi \in \GZK$. On the other hand, it is readily checked that $C_0$,
the critical submanifold of $S$ with critical value $0$, has homology
$$H_{\Z_k, p}(C_0) = H_{\Z_k, p}(pt) = H_{p}(B\Z_k) = \Z_k$$
for all $p \geq 0$ and any $\phi \in \GZK$. This is
because the $\Z_k$-action in this example fixes the point at $\infty$ (South pole in the base $S^{2n}$) and
by perturbing the generating function  $S \geq 0$
for $\phi \in \G$ one can assume a single critical point at $\infty$ with critical value
$0$. Thus, as Sandon computes for certain special $\phi$ supported in $\widehat{B}(R)$
\begin{equation}
G_{\Z_k, p}^{(\epsilon_\phi, \infty]}(\phi)  = \Z_k \label{eq:equivariant}
\end{equation}
for all $p \geq 2n$, when 
$\epsilon_\phi$ is sufficiently small that a generic\footnote{(see footnote 5  - all references to the surrogate in Sections 3 and 4 were updated accordingly Feb. 2015)} generating function $S$ for $\phi$ has
no nonzero critical values less than or equal to $\epsilon_\phi$.
In fact in this setting, for all compactly supported $\phi$, by the equivariant Thom isomorphism,
$G_{\Z_k, p}^{(-\infty, \infty]}(\phi)$ is given by 
$H_{\Z_k, *}(S^{2n})$ and so has a single generator in degrees $p = 0, \ldots, 2n-1$
and two generators in degrees $p \geq 2n$. Computation with the long exact sequence for 
the triple $(E, E^{\epsilon}, E^{-1})$ shows that
a single generator $\kappa_p$ in all degrees $p \geq 0$ is lost when passing action threshold 0,
hence $c(\phi, \kappa_p) = 0$ for all $p \geq 0$; denote their span $P := \langle \kappa_p: p \geq 0 \rangle \subset G_{\Z_k, p}^{(-\infty, \infty]}(\phi)$. Now define
the surrogate 
$G_{\Z_k, p}^{(0_+, \infty]}(\phi) := G_{\Z_k, p}^{(-\epsilon, \infty]}(\phi)/P$.
By construction this vector space is conjugation-invariant. 
Then 
for fixed positive integer $p$ put
$$W_{a} (\phi) := \left \{
\begin{array}{ll}
G_{\Z_k, p}^{(a, \infty]}(\phi) &\text{ if } a < 0 \\ 
G_{\Z_k, p}^{(0_+, \infty]}(\phi) &\text{ if } a=0  \\
G_{\Z_k, p}^{(a, \infty]}(\phi) &\text{ if } a>0.
\end{array}
\right .
$$
This is a persistence module
with persistence maps $(w_a^b)_*: W_{a}(\phi) \to W_{b}(\phi)$ for $a < b$,
defined as before. 
Recall that $W_0$ has a single generator in each degree $p \geq 2n$;
denote these $\eta^p$ and define
$$c^p_{\Z_k}(\phi) :=  \rho(\eta^p),$$
the persistence of $\eta^p$. 
The conjugation-invariance
of the groups $W_a(\phi)$ for $a \in \Z$ implies
$[c^p_{\Z_k}(\phi)]$ is conjugation-invariant. By taking the supremum
over all $\phi$ supported in a bounded domain $\mathcal W \subset \R^{2n} \times S^1$ one obtains a contact 
invariant $[c^p_{\Z_k}(\mathcal W)]$ which is by \mbox{definition monotone}. From the results in 
Sandon \cite{sandon-equiv} one computes
\begin{equation}
[c^p_{\Z_k}(\widehat {B(R)}] = [\ell R] \text{ when }p = 2n\ell \label{eq-comp}
\end{equation}
so $[c^p_{\Z_k}]$ is non-trivial and may be viewed as a {\bf $\Z_k$-equivariant contact capacity}.
Moreover, this capacity is sufficient to establish that squeezing via (compactly supported) $\Z_k$-equivariant contactomorphisms
is impossible in $\R^{2n} \times S^1$ and thus to {\bf prove orderability of lens spaces} 
which Milin \cite{milin} and Sandon \cite{sandon-equiv} proved 
by means of contact homology 
and equivariant GF homology groups respectively. 

\bigskip
\paragraph{\sc Remark:}
The interest of $[c^p_{\Z_k}]$ 
is primarily in confirming that 
the simpler capacity-based non-squeezing argument 
of the non-equivariant setting is also possible in the equivariant one
and to view such invariants more generally from a persistence viewpoint.
This argument is only simpler in a minor way,
as it makes use of the existing framework of Sandon \cite{sandon-equiv}:
the existence and uniqueness of $\Z_k$-equivariant generating functions 
for this setting (which she establishes by 
a $\Z_k$-equivariant adaptation of Chaperon and Th\'eret's arguments
and uses 
to show her $\Z_k$-equivariant GF homology groups are well-defined)
and also the functorial properties of these groups 
which in our case imply they form a $\Z_k$-invariant persistence module.
Regarding Sandon's computation of the groups $W_a, a \in \R$ for each $p \in \N$
which we use in \eqref{eq-comp} there is primarily one difficult case,
namely when $a$ passes critical points of GF index $2n\ell$ and $2n(\ell - 1)$.
She handles that case by perturbing the generating function to obtain a true Morse function
and arguing directly with the Morse complex. A similar technique is also used by
Milin. 

\bigskip
\paragraph{\sc II) Symplectic version} 
In fact Sandon also defines $\Z_k$-equivariant generating functions for Hamiltonian diffeomorphisms
of $R^{2n}$ for the same $\Z_k$-action. As above we may define 
symplectic capacities $c^p_{\Z_k}(\cdot), p \in \N$. Unlike the contact version just defined,
these do not
give anything beyond the Viterbo capacity $c_V$ at least for balls
(because $c_V$ is a real-valued symplectic invariant which already separates all balls). However, we do have
$$c^p_{\Z_k, \text{contact}}(\widehat {\mathcal U}) =c^p_{\Z_k, \text{symp}}(\mathcal U)$$
for $\widehat {\mathcal U} := \mathcal U \times S^1$, $\mathcal U \subset \R^{2n}$ 
which was used in computing \eqref{eq-comp}.

\paragraph{Acknowledgements} These informal notes were 
circulated to a small group of people in January 2014.
Due to recent interest I am making them available on arXiv.
They were written in the months
after completing my PhD in Computer Science (University of Chicago, August 2013).
A portion of my PhD thesis was in topological data analysis and
my initial advisers were Partha Niyogi and Steve Smale. I thank them both
as well as Shmuel Weinberger, Vin de Silva and Fred Chazal 
for useful discussions concerning persistent
homology. At the same time, during my studies at Chicago I took a
reading course from Leonid Polterovich in the Mathematics Department
focusing on recent work of Sheila Sandon and I thank him for
his mathematical guidance and for many interesting discussions. 
Having one foot in the field of
topological data analysis and the other in the field of symplectic/contact geometry 
gave rise to the observations in these notes and the
work of Chazal et al \cite{vin} was a key abstraction which made it possible
to formulate the connection. A forthcoming paper \cite{fraser-new}
gives a concrete application of this viewpoint.

\end{document}